\documentclass[11pt,reqno]{amsart}
\usepackage{mathrsfs}
\usepackage{}

\def\subjclass#1{{\renewcommand{\thefootnote}{}%
\footnote{\emph{Mathematics Subject Classification (2020):} #1}}}

\usepackage{amsfonts}
\usepackage{amssymb}

\DeclareMathOperator{\divg}{div}

\date{\today}

\hoffset=- 2cm \voffset=- 0cm 
 \theoremstyle{plain}
\newtheorem{Thm}{Theorem}

\newtheorem{Rem}[Thm]{Remark}
\newtheorem{Lem}[Thm]{Lemma}
\newtheorem{Cor}[Thm]{Corollary}

\usepackage{amssymb}
\usepackage{color}

\def\0{\mathbf 0}

\def\v{\vskip}

\errorcontextlines=0 \numberwithin{equation}{section}
\numberwithin{Thm}{section}
\allowdisplaybreaks

\begin{document}
\large

\title[Liouville type theorems]
{Some remarks on Liouville type theorems for the 3D steady tropical climate model}

\author{Yanyan Dong}

\address{Yanyan Dong: School of Mathematics and Physics, Anqing Normal University, Anqing 246133, China}
\email{2981165080@qq.com}

\author{Zhibing Zhang}

\address{Zhibing Zhang: School of Mathematics and Physics, Key Laboratory of Modeling, Simulation and Control of Complex Ecosystem in Dabie Mountains of Anhui Higher Education Institutes, Anqing Normal University, Anqing 246133, China}
\email{zhibingzhang29@126.com}%

\thanks{}

\keywords{Liouville type theorems, tropical climate model, Poincar{\'{e}}-Sobolev inequality}

\subjclass{35B53, 35Q35, 35A02}

\begin{abstract}
Observing the special structure of the system and using the Poincar{\'{e}}-Sobolev inequality, we establish Liouville type theorems for  the 3D steady tropical climate model under certain conditions on $u$, $v$, $\nabla \theta$. Our results extend and improve a Liouville type result of Cho-In-Yang (arXiv:2312.17441). 
\end{abstract}
\maketitle

\section{Introduction}
We consider the following stationary tropical climate model in three dimension
\begin{equation}\label{equ1.1}
  \left\{
    \begin{array}{ll}
     -\Delta u+(u\cdot\nabla) u+\nabla \pi+\divg(v\otimes v)=0,  \\
   -\Delta v+(u\cdot\nabla) v+\nabla \theta+(v\cdot\nabla)u=0, \\
   -\Delta \theta+u\cdot\nabla\theta+\divg v=0,\\
     \divg u=0,
\end{array}
  \right.
\end{equation}
where $u= (u_1 , u_2 , u_3 )$ is the barotropic mode, $v= (v_1 , v_2 , v_3 )$ is the first baroclinic mode of vector velocity, $\theta$ is the temperature, and $\pi$ is the pressure.

To our knowledge, there are relatively few Liouville type results for the tropical climate models so far. Ding and Wu \cite{FW21} investigated the tropical climate model and proved Liouville type theorems under one of the following assumptions
$$
\aligned
&(\mathrm{i})~~~~u\in L^p(\mathbb{R}^{3}) ~~~~~~~~~~~~\text{and}~~~~~~~~~~~~v,\theta\in L^p(\mathbb{R}^3)\cap L^3(\mathbb{R}^{3}) ~~~~~~~~~~~~\text{with}~~~~~~~~~~~~ 3\leq p\leq \frac{9}{2},\\
&(\mathrm{ii})~~~~u\in L^{p,q}(\mathbb{R}^3) ~~~~~~~~~~~~\text{and}~~~~~~~~~~~~ v, \theta \in L^{p,q}(\mathbb{R}^3)\cap L^{3,q} ~~~~~~~~~~~~\text{with}~~~~~~~~~~~~ 3 \leq p\leq\frac{9}{2},\;3<q<+\infty,\\
&(\mathrm{iii})~~~~u\in L^{3}(\mathbb{R}^3), v \in L^{2}(\mathbb{R}^3)~~~~~~~~~~~~\text{and}~~~~~~~~~~~~\nabla u,\nabla v,\nabla\theta \in L^{2}(\mathbb{R}^3),
\endaligned
$$
where $L^{p,q}(\mathbb{R}^3)$ is the Lorentz space. Later, Yuan and Wang \cite{BY23} extended their results to local Morrey spaces. Cho et al. \cite{CIY23} improved the result of Ding and Wu by removing the condition $\nabla u,\nabla v \in L^{2}(\mathbb{R}^3)$ from the above assumption $(\mathrm{iii})$.
Recently, Cho et al. \cite{CIY24} established an improved Liouville type theorem for \eqref{equ1.1} under one of the following conditions
\begin{subequations}
\begin{align}
&(\mathrm{i})\;p\in\left(\frac{3}{2},3\right),\;q,r\in[1,2),\;\liminf_{R\rightarrow+\infty}(G(u;p,R)+H(v;q,R)+H(\theta;r,R))<+\infty, \label{1.2a}\\
&(\mathrm{ii})\liminf_{R\rightarrow+\infty}(G(u;3,R)+H(v;2,R))=0,\;\limsup_{R\rightarrow+\infty}H(\theta;2,R)<+\infty,\label{1.2b}\\
&(\mathrm{iii})\liminf_{R\rightarrow+\infty}(G(u;3,R)+H(\theta;2,R))=0,\;\limsup_{R\rightarrow+\infty}H(v;2,R)<+\infty,\label{1.2c}
\end{align}
\end{subequations}
where the functions $G$ and $H$ are defined by
$$ G(f;p,R):=R^{-\left(\frac{2}{p}-\frac{1}{3}\right)}\|f\|_{L^p(B_{2R}\backslash B_{R/2})}
~~~~~~~~~~~~and~~~~~~~~~~~~
H(f;p,R):=R^{-\left(\frac{3}{2p}-\frac{1}{4}\right)}\|f\|_{L^p(B_{2R}\backslash B_{R/2})}.
$$
Very recently, Fang and Zhang \cite{FZ25} extended and improved the work of Cho et al. \cite{CIY24} by combining an iteration argument and some delicate estimates of several integrals.
In order to show Liouville type results conveniently, we introduce some notations. For simplicity, we denote
$$
\aligned
&X_{\alpha,p,R}=R^{-\alpha}\|u\|_{L^p\left(B_{2R} \backslash B_R\right)},\;Y_{\beta,q,R}=R^{-\beta}\|v\|_{L^q\left(B_{2R} \backslash B_R\right)},\;Z_{\gamma,r,R}=R^{-\gamma}\|\theta\|_{L^r\left(B_{2R} \backslash B_R\right)},\\
&\overline{Z}_{\gamma,r,R}=R^{-\gamma}\|\theta-\overline{\theta}_{R}\|_{L^r\left(B_{2R} \backslash B_R\right)},\;Z'_{\gamma,\lambda,R}=R^{-\gamma}\|\nabla\theta\|_{L^\lambda\left(B_{2R} \backslash B_R\right)},
\endaligned
$$
where $\overline{\theta}_{R}$ is the mean value of $\theta$ over $B_{2R}\backslash B_R$. For the assumptions $\mathrm{(A1)}$-$\mathrm{(A17)}$ and $\mathrm{(B1)}$-$\mathrm{(B10)}$, $\mathrm{(A1)'}$-$\mathrm{(A17)'}$ and  $\mathrm{(B1)'}$-$\mathrm{(B10)'}$, see Section \ref{sec2}.
More precisely, Fang and Zhang \cite{FZ25} obtained the following Liouville type theorems:

\begin{itemize}
\item[(i)] Suppose $p\in[3,\frac{9}{2}]$, $q,r\in [1,6]$ and $\alpha\in\left[0,\frac{3}{p}-\frac{2}{3}\right]$, $\beta\in\left[0,\frac{3}{q}-\frac{1}{2}\right]$, $\gamma\in\left[0,\frac{3}{r}-\frac{1}{2}\right]$. Assume that there exists a sequence $R_j\nearrow+\infty$ such that
\begin{equation*}
\lim\limits_{j\rightarrow+\infty}X_{\alpha,p,R_j}=0,\;\lim\limits_{j\rightarrow+\infty}Y_{\beta,q,R_j}<+\infty,\;\lim\limits_{j\rightarrow+\infty}Z_{\gamma,r,R_j}<+\infty.
\end{equation*}
Moreover, if $p,q,r,\alpha,\beta,\gamma$ satisfy one of the nine assumptions $\mathrm{(A1)}$, $\mathrm{(A2)}$, $\cdots$, $\mathrm{(A9)}$, then $u=v=0$ and $\theta=0$.
\item[(ii)] Suppose $p\in(\frac{3}{2},3)$, $q,r\in [1,6]$ and $\alpha\in\left[0,\frac{2}{p}-\frac{1}{3}\right]$, $\beta\in\left[0,\frac{3}{q}-\frac{1}{2}\right]$, $\gamma\in\left[0,\frac{3}{r}-\frac{1}{2}\right]$. Assume that
$$\liminf\limits_{R\rightarrow+\infty}\left(X_{\alpha,p,R}+Y_{\beta,q,R}+Z_{\gamma,r,R}\right)<+\infty.$$
Moreover, if $p,q,r,\alpha,\beta,\gamma$ satisfy one of the eight assumptions $\mathrm{(A10)}$, $\mathrm{(A11)}$, $\cdots$, $\mathrm{(A17)}$, then $u=v=0$ and $\theta=0$.
\end{itemize}
As a corollary, Liouville type theorem holds if $u\in L^p(\mathbb{R}^3)$, $v\in L^q(\mathbb{R}^3)$, $\theta\in L^r(\mathbb{R}^3)$, where $p,q,r$ satisfy one of the ten assumptions $\mathrm{(B1)}$, $\mathrm{(B2)}$, $\cdots$, $\mathrm{(B10)}$.

Our goal of this paper is to extend and improve the result of Cho et al. \cite{CIY23} and to establish more Liouville type theorems under the conditions on $u$, $v$, $\nabla \theta$. Observing the special structure of the system for $\theta$ and using the Poincar{\'{e}}-Sobolev inequality, we achieve this goal. Actually, we obtain Liouville type results under weaker conditions proposed on $u$, $v$, $\theta-\overline{\theta}_{R}$. Besides, the result of Cho et al. \cite{CIY23} can be viewed as a special case of Corollary \ref{Cor1.5} with $p,q,\lambda$ satisfying $\mathrm{(B6)'}$.

Our main results are stated as follows.
\begin{Thm}\label{main1}
Let $(u,\pi,v,\theta)$ be a smooth solution of \eqref{equ1.1}. Suppose $p\in[3,\frac{9}{2}]$, $q,r\in [1,6]$ and $\alpha\in\left[0,\frac{3}{p}-\frac{2}{3}\right]$, $\beta\in\left[0,\frac{3}{q}-\frac{1}{2}\right]$, $\gamma\in\left[0,\frac{3}{r}-\frac{1}{2}\right]$. Assume that there exists a sequence $R_j\nearrow+\infty$ such that
\begin{equation}\label{lim}
\lim\limits_{j\rightarrow+\infty}X_{\alpha,p,R_j}=0,\;\lim\limits_{j\rightarrow+\infty}Y_{\beta,q,R_j}<+\infty,\;\lim\limits_{j\rightarrow+\infty}\overline{Z}_{\gamma,r,R_j}<+\infty.
\end{equation}
Moreover, if $p,q,r,\alpha,\beta,\gamma$ satisfy one of the nine assumptions $\mathrm{(A1)}$, $\mathrm{(A2)}$, $\cdots$, $\mathrm{(A9)}$, then $u=v=0$ and $\theta$ is a constant.
\end{Thm}

\begin{Thm}\label{main2}
Let $(u,\pi,v,\theta)$ be a smooth solution of \eqref{equ1.1}. Suppose $p\in(\frac{3}{2},3)$, $q,r\in [1,6]$ and $\alpha\in\left[0,\frac{2}{p}-\frac{1}{3}\right]$, $\beta\in\left[0,\frac{3}{q}-\frac{1}{2}\right]$, $\gamma\in\left[0,\frac{3}{r}-\frac{1}{2}\right]$. Assume that
$$\liminf\limits_{R\rightarrow+\infty}\left(X_{\alpha,p,R}+Y_{\beta,q,R}+\overline{Z}_{\gamma,r,R}\right)<+\infty.$$
Moreover, if $p,q,r,\alpha,\beta,\gamma$ satisfy one of the eight assumptions $\mathrm{(A10)}$, $\mathrm{(A11)}$, $\cdots$, $\mathrm{(A17)}$, then $u=v=0$ and $\theta$ is a constant.
\end{Thm}

As consequences of Theorem \ref{main1}, Theorem \ref{main2} and the Poincar{\'{e}}-Sobolev inequality (see Lemma \ref{Lem2.2})
$$\|\theta-\overline{\theta}_{R}\|_{L^r(B_{2R}\backslash B_R)}\leq C\|\nabla\theta\|_{L^\lambda (B_{2R}\backslash B_R)} \text{ with }\frac{1}{r}=\frac{1}{\lambda}-\frac{1}{3},\;1\leq\lambda<3,$$
we have the following two corollaries.

\begin{Cor}\label{Cor1.3}
Let $(u,\pi,v,\theta)$ be a smooth solution of \eqref{equ1.1}.
Suppose $p\in[3,\frac{9}{2}]$, $q\in [1,6]$, $\lambda\in [1,2]$ and $\alpha\in\left[0,\frac{3}{p}-\frac{2}{3}\right]$, $\beta\in\left[0,\frac{3}{q}-\frac{1}{2}\right]$, $\gamma\in\left[0,\frac{3}{\lambda}-\frac{3}{2}\right]$. Assume that there exists a sequence $R_j\nearrow+\infty$ such that
\begin{equation*}
\lim\limits_{j\rightarrow+\infty}X_{\alpha,p,R_j}=0,\;\lim\limits_{j\rightarrow+\infty}Y_{\beta,q,R_j}<+\infty,\;\lim\limits_{j\rightarrow+\infty}Z'_{\gamma,\lambda,R_j}<+\infty.
\end{equation*}
Moreover, if $p,q,\lambda,\alpha,\beta,\gamma$ satisfy one of the nine assumptions $\mathrm{(A1)'}$, $\mathrm{(A2)'}$, $\cdots$, $\mathrm{(A9)'}$, then $u=v=0$ and $\theta$ is a constant.
\end{Cor}

\begin{Cor}\label{Cor1.4}
Let $(u,\pi,v,\theta)$ be a smooth solution of \eqref{equ1.1}.
Suppose $p\in(\frac{3}{2},3)$, $q\in [1,6]$, $\lambda\in [1,2]$ and $\alpha\in\left[0,\frac{2}{p}-\frac{1}{3}\right]$, $\beta\in\left[0,\frac{3}{q}-\frac{1}{2}\right]$, $\gamma\in\left[0,\frac{3}{\lambda}-\frac{3}{2}\right]$. Assume that
$$\liminf\limits_{R\rightarrow+\infty}\left(X_{\alpha,p,R}+Y_{\beta,q,R}+Z'_{\gamma,\lambda,R}\right)<+\infty.$$
Moreover, if $p,q,\lambda,\alpha,\beta,\gamma$ satisfy one of the eight assumptions $\mathrm{(A10)'}$, $\mathrm{(A11)'}$, $\cdots$, $\mathrm{(A17)'}$, then $u=v=0$ and $\theta$ is a constant.
\end{Cor}

Particularly, taking $\alpha=\beta=\gamma=0$ in Corollary \ref{Cor1.3} and Corollary \ref{Cor1.4}, and noticing some endpoint cases in Remark \ref{Rem3.1}, after a careful calculation and classification, we have the following conclusion.
\begin{Cor}\label{Cor1.5}
	Let $(u,\pi,v,\theta)$ be a smooth solution of \eqref{equ1.1}. Suppose $u\in L^p(\mathbb{R}^3)$, $v\in L^q(\mathbb{R}^3)$, $\nabla\theta\in L^\lambda(\mathbb{R}^3)$, where $p,q,\lambda$ satisfy one of the ten assumptions $\mathrm{(B1)'}$, $\mathrm{(B2)'}$, $\cdots$, $\mathrm{(B10)'}$. Then $u=v=0$ and $\theta$ is a constant.
\end{Cor}

Finally, we give an arrangement of the rest of this paper. In Section \ref{sec2}, we introduce some basic assumptions on the parameters, the property of the Bogovskii operator and some key lemmas. Section \ref{sec3} is devoted to proving Liouville type theorems.
\v0.1in
\section{Preliminaries}\label{sec2}

Throughout this article, we use $C$ to denote a finite inessential constant which may be different from line to line.
Let $p'$ denote the conjugate exponent to $p$, i.e., $p'=\frac{p}{p-1}$. When we consider the quantities $X_{\alpha,p,R}$, $Y_{\beta,q,R}$, $Z_{\gamma,r,R}$  or $X_{\alpha,p,R}$, $Y_{\beta,q,R}$, $\overline{Z}_{\gamma,r,R}$, the parameters $p,q,r,\alpha,\beta,\gamma$ satisfy one of the following assumptions
\begin{align*}
\mathrm{(A1)}&\;p\in\left[3,\frac{9}{2}\right],q\in[1,2),r\in[1,2),\\
&\alpha+\frac{(4p-6)q}{(6-q)p}\beta\leq1,\;\alpha+\frac{(4p-6)r}{(6-r)p}\gamma\leq1,\;\frac{2q}{6-q}\beta+\frac{2r}{6-r}\gamma\leq1;\\
\mathrm{(A2)}&\;p\in\left[3,\frac{9}{2}\right],q\in[1,2),r\in[2,2p'),\\
&\alpha+\frac{(4p-6)q}{(6-q)p}\beta\leq1,\;\alpha+\frac{(4p-6)r}{(6-r)p}\gamma\leq1,\;\frac{2q}{6-q}\beta+\gamma\leq\frac{6-r}{2r};\\
\mathrm{(A3)}&\;p\in\left[3,\frac{9}{2}\right],q\in[1,2),r\in[2p',6],\\
&\alpha+\frac{(4p-6)q}{(6-q)p}\beta\leq1,\;\alpha+2\gamma\leq\frac{3}{p}+\frac{6}{r}-2,\;\frac{2q}{6-q}\beta+\gamma\leq\frac{6-r}{2r};\\
\mathrm{(A4)}&\;p\in\left[3,\frac{9}{2}\right],q\in[2,2p'),r\in[1,2),\\
&\alpha+\frac{(4p-6)q}{(6-q)p}\beta\leq1,\;\alpha+\frac{(4p-6)r}{(6-r)p}\gamma\leq1,\;\beta+\frac{2r}{6-r}\gamma\leq\frac{6-q}{2q};\\
\mathrm{(A5)}&\;p\in\left[3,\frac{9}{2}\right],q\in[2,2p'),r\in[2,2p'),\\
&\alpha+\frac{(4p-6)q}{(6-q)p}\beta\leq1,\;\alpha+\frac{(4p-6)r}{(6-r)p}\gamma\leq1,\;\beta+\gamma<\frac{3}{q}+\frac{3}{r}-2;\\
\mathrm{(A6)}&\;p\in\left[3,\frac{9}{2}\right],q\in[2,2p'),r\in[2p',6],\\
&\alpha+\frac{(4p-6)q}{(6-q)p}\beta\leq1,\;\alpha+2\gamma\leq\frac{3}{p}+\frac{6}{r}-2,\;\beta+\gamma<\frac{3}{q}+\frac{3}{r}-2;\\
\mathrm{(A7)}&\;p\in\left[3,\frac{9}{2}\right],q\in[2p',6],r\in[1,2),\\
&\alpha+2\beta\leq\frac{3}{p}+\frac{6}{q}-2,\;\alpha+\frac{(4p-6)r}{(6-r)p}\gamma\leq1,\;\beta+\frac{2r}{6-r}\gamma\leq\frac{6-q}{2q};\\
\mathrm{(A8)}&\;p\in\left[3,\frac{9}{2}\right],q\in[2p',6],r\in[2,2p'),\\
&\alpha+2\beta\leq\frac{3}{p}+\frac{6}{q}-2,\;\alpha+\frac{(4p-6)r}{(6-r)p}\gamma\leq1,\;\beta+\gamma<\frac{3}{q}+\frac{3}{r}-2;\\
\mathrm{(A9)}&\;p\in\left[3,\frac{9}{2}\right],q\in[2p',6],r\in[2p',6],\\
&\alpha+2\beta\leq\frac{3}{p}+\frac{6}{q}-2,\;\alpha+2\gamma\leq\frac{3}{p}+\frac{6}{r}-2,\;\beta+\gamma<\frac{3}{q}+\frac{3}{r}-2;\\
\mathrm{(A10)}&\;p\in\left(\frac{3}{2},3\right),q\in[1,2),r\in[1,2),\\
&\alpha+\frac{(4p-6)q}{(6-q)p}\beta\leq1,\;\alpha+\frac{(4p-6)r}{(6-r)p}\gamma\leq1,\;\frac{2q}{6-q}\beta+\frac{2r}{6-r}\gamma\leq1;\\
\mathrm{(A11)}&\;p\in\left(\frac{3}{2},3\right),q\in[1,2),r\in[2,2p'),\\
&\alpha+\frac{(4p-6)q}{(6-q)p}\beta\leq1,\;\alpha+\frac{(4p-6)r}{(6-r)p}\gamma\leq1,\;\frac{2q}{6-q}\beta+\gamma\leq\frac{6-r}{2r};\\
\mathrm{(A12)}&\;p\in\left(\frac{3}{2},3\right),q\in[1,2),r\in[2p',6],\\
&\alpha+\frac{(4p-6)q}{(6-q)p}\beta\leq1,\;\alpha+2\gamma<\frac{3}{p}+\frac{6}{r}-2,\;\frac{2q}{6-q}\beta+\gamma\leq\frac{6-r}{2r};\\
\mathrm{(A13)}&\;p\in\left(\frac{3}{2},3\right),q\in[2,2p'),r\in[1,2),\\
&\alpha+\frac{(4p-6)q}{(6-q)p}\beta\leq1,\;\alpha+\frac{(4p-6)r}{(6-r)p}\gamma\leq1,\;\beta+\frac{2r}{6-r}\gamma\leq\frac{6-q}{2q};\\
\mathrm{(A14)}&\;p\in\left(\frac{3}{2},3\right),q\in[2,2p'),r\in[2,2p'),\\
&\alpha+\frac{(4p-6)q}{(6-q)p}\beta\leq1,\;\alpha+\frac{(4p-6)r}{(6-r)p}\gamma\leq1,\;\beta+\gamma<\frac{3}{q}+\frac{3}{r}-2;\\
\mathrm{(A15)}&\;p\in\left(\frac{3}{2},3\right),q\in[2,2p'),r\in[2p',6],\\
&\alpha+\frac{(4p-6)q}{(6-q)p}\beta\leq1,\;\alpha+2\gamma<\frac{3}{p}+\frac{6}{r}-2,\;\beta+\gamma<\frac{3}{q}+\frac{3}{r}-2;\\
\mathrm{(A16)}&\;p\in\left(\frac{3}{2},3\right),q\in[2p',6],r\in[1,2),\\
&\alpha+2\beta<\frac{3}{p}+\frac{6}{q}-2,\;\alpha+\frac{(4p-6)r}{(6-r)p}\gamma\leq1,\;\beta+\frac{2r}{6-r}\gamma\leq\frac{6-q}{2q};\\
\mathrm{(A17)}&\;p\in\left(\frac{3}{2},3\right),q\in[2p',6],r\in[2,2p'),\\
&\alpha+2\beta<\frac{3}{p}+\frac{6}{q}-2,\;\alpha+\frac{(4p-6)r}{(6-r)p}\gamma\leq1,\;\beta+\gamma<\frac{3}{q}+\frac{3}{r}-2.
\end{align*}
Especially, when $\alpha=\beta=\gamma=0$, $p,q,r$ satisfy one of the following assumptions
\begin{align*}
\mathrm{(B1)}&\;p\in\left[3,\frac{9}{2}\right],q\in[1,2p'),r\in[1,2p');\\
\mathrm{(B2)}&\;p\in\left(\frac{3}{2},3\right),q\in[1,2),r\in[1,6];\\
\mathrm{(B3)}&\;p\in\left(\frac{3}{2},3\right),q\in[1,6],r\in[1,2);\\
\mathrm{(B4)}&\;p\in\left[3,\frac{9}{2}\right],q\in[1,2),r\in[2p',6],\text{ with }\;\frac{1}{p}+\frac{2}{r}\geq\frac{2}{3};\\
\mathrm{(B5)}&\;p\in\left[3,\frac{9}{2}\right],q\in[2p',6],r\in[1,2),\text{ with }\;\frac{1}{p}+\frac{2}{q}\geq\frac{2}{3};\\
\mathrm{(B6)}&\;p\in\left[3,\frac{9}{2}\right],q\in[2,2p'),r\in[2p',6],\text{ with }\;\frac{1}{p}+\frac{2}{r}\geq\frac{2}{3},\;\frac{1}{q}+\frac{1}{r}\geq\frac{2}{3};\\
\mathrm{(B7)}&\;p\in\left[3,\frac{9}{2}\right],q\in[2p',6],r\in[2,2p'),\text{ with }\;\frac{1}{p}+\frac{2}{q}\geq\frac{2}{3},\;\frac{1}{q}+\frac{1}{r}\geq\frac{2}{3};\\
\mathrm{(B8)}&\;p\in\left[3,\frac{9}{2}\right],q\in[2p',6],r\in[2p',6],\text{ with }\; \frac{1}{p}+\frac{2}{q}\geq\frac{2}{3},\;\frac{1}{p}+\frac{2}{r}\geq\frac{2}{3},\;\frac{1}{q}+\frac{1}{r}\geq\frac{2}{3};\\
\mathrm{(B9)}&\;p\in\left(\frac{3}{2},3\right),q\in[2,2p'),r\in[2,6],\text{ with }\;\frac{1}{q}+\frac{1}{r}\geq\frac{2}{3};\\
\mathrm{(B10)}&\;p\in\left(\frac{3}{2},3\right),q\in[2,6],r\in[2,2p'),\text{ with }\;\frac{1}{q}+\frac{1}{r}\geq\frac{2}{3}.\\
\end{align*}

When we consider the quantities $X_{\alpha,p,R}$, $Y_{\beta,q,R}$, $Z'_{\gamma,r,R}$, the parameters $p$, $q$, $\lambda$, $\alpha$, $\beta$, $\gamma$ satisfy one of the following assumptions
\begin{align*}
\mathrm{(A1)'}&\;p\in\left[3,\frac{9}{2}\right],\;q\in[1,2),\;\lambda\in\left[1,\frac{6}{5}\right),\\
&\alpha+\frac{(4p-6)q}{(6-q)p}\beta\leq1,\;\alpha+\frac{(4p-6)\lambda}{(6-3\lambda)p}\gamma\leq1,\;\frac{2q}{6-q}\beta+\frac{2\lambda}{6-3\lambda}\gamma\leq1;\\
\mathrm{(A2)'}&\;p\in\left[3,\frac{9}{2}\right],\;q\in[1,2),\;\lambda\in\left[\frac{6}{5},\frac{6p}{5p-3}\right),\\
&\alpha+\frac{(4p-6)q}{(6-q)p}\beta\leq1,\;\alpha+\frac{(4p-6)\lambda}{(6-3\lambda)p}\gamma\leq1,\;\frac{2q}{6-q}\beta+\gamma\leq\frac{6-3\lambda}{2\lambda};\\
\mathrm{(A3)'}&\;p\in\left[3,\frac{9}{2}\right],\;q\in[1,2),\;\lambda\in\left[\frac{6p}{5p-3},2\right],\\
&\alpha+\frac{(4p-6)q}{(6-q)p}\beta\leq1,\;\alpha+2\gamma\leq\frac{3}{p}+\frac{6}{\lambda}-4,\;\frac{2q}{6-q}\beta+\gamma\leq\frac{6-3\lambda}{2\lambda};\\
\mathrm{(A4)'}&\;p\in\left[3,\frac{9}{2}\right],\;q\in[2,2p'),\;\lambda\in\left[1,\frac{6}{5}\right),\\
&\alpha+\frac{(4p-6)q}{(6-q)p}\beta\leq1,\;\alpha+\frac{(4p-6)\lambda}{(6-3\lambda)p}\gamma\leq1,\;\beta+\frac{2\lambda}{6-3\lambda}\gamma\leq\frac{6-q}{2q};\\
\mathrm{(A5)'}&\;p\in\left[3,\frac{9}{2}\right],\;q\in[2,2p'),\;\lambda\in\left[\frac{6}{5},\frac{6p}{5p-3}\right),\\
&\alpha+\frac{(4p-6)q}{(6-q)p}\beta\leq1,\;\alpha+\frac{(4p-6)\lambda}{(6-3\lambda)p}\gamma\leq1,\;\beta+\gamma<\frac{3}{q}+\frac{3}{\lambda}-3;\\
\mathrm{(A6)'}&\;p\in\left[3,\frac{9}{2}\right],\;q\in[2,2p'),\;\lambda\in\left[\frac{6p}{5p-3},2\right],\\
&\alpha+\frac{(4p-6)q}{(6-q)p}\beta\leq1,\;\alpha+2\gamma\leq\frac{3}{p}+\frac{6}{\lambda}-4,\;\beta+\gamma<\frac{3}{q}+\frac{3}{\lambda}-3;\\
\mathrm{(A7)'}&\;p\in\left[3,\frac{9}{2}\right],\;q\in[2p',6],\;\lambda\in\left[1,\frac{6}{5}\right),\\
&\alpha+2\beta\leq\frac{3}{p}+\frac{6}{q}-2,\;\alpha+\frac{(4p-6)\lambda}{(6-3\lambda)p}\gamma\leq1,\;\beta+\frac{2\lambda}{6-3\lambda}\gamma\leq\frac{6-q}{2q};\\
\mathrm{(A8)'}&\;p\in\left[3,\frac{9}{2}\right],\;q\in[2p',6],\;\lambda\in\left[\frac{6}{5},\frac{6p}{5p-3}\right),\\
&\alpha+2\beta\leq\frac{3}{p}+\frac{6}{q}-2,\;\alpha+\frac{(4p-6)\lambda}{(6-3\lambda)p}\gamma\leq1,\;\beta+\gamma<\frac{3}{q}+\frac{3}{\lambda}-3;\\
\mathrm{(A9)'}&\;p\in\left[3,\frac{9}{2}\right],\;q\in[2p',6],\;\lambda\in\left[\frac{6p}{5p-3},2\right],\\
&\alpha+2\beta\leq\frac{3}{p}+\frac{6}{q}-2,\;\alpha+2\gamma\leq\frac{3}{p}+\frac{6}{\lambda}-4,\;\beta+\gamma<\frac{3}{q}+\frac{3}{\lambda}-3;\\
\mathrm{(A10)'}&\;p\in\left(\frac{3}{2},3\right),\;q\in[1,2),\;\lambda\in\left[1,\frac{6}{5}\right),\\
&\alpha+\frac{(4p-6)q}{(6-q)p}\beta\leq1,\;\alpha+\frac{(4p-6)\lambda}{(6-3\lambda)p}\gamma\leq1,\;\frac{2q}{6-q}\beta+\frac{2\lambda}{6-3\lambda}\gamma\leq1;\\
\mathrm{(A11)'}&\;p\in\left(\frac{3}{2},3\right),\;q\in[1,2),\;\lambda\in\left[\frac{6}{5},\frac{6p}{5p-3}\right),\\
&\alpha+\frac{(4p-6)q}{(6-q)p}\beta\leq1,\;\alpha+\frac{(4p-6)\lambda}{(6-3\lambda)p}\gamma\leq1,\;\frac{2q}{6-q}\beta+\gamma\leq\frac{6-3\lambda}{2\lambda};\\
\mathrm{(A12)'}&\;p\in\left(\frac{3}{2},3\right),\;q\in[1,2),\;\lambda\in\left[\frac{6p}{5p-3},2\right],\\
&\alpha+\frac{(4p-6)q}{(6-q)p}\beta\leq1,\;\alpha+2\gamma<\frac{3}{p}+\frac{6}{\lambda}-4,\;\frac{2q}{6-q}\beta+\gamma\leq\frac{6-3\lambda}{2\lambda};\\
\mathrm{(A13)'}&\;p\in\left(\frac{3}{2},3\right),\;q\in[2,2p'),\;\lambda\in\left[1,\frac{6}{5}\right),\\	&\alpha+\frac{(4p-6)q}{(6-q)p}\beta\leq1,\;\alpha+\frac{(4p-6)\lambda}{(6-3\lambda)p}\gamma\leq1,\;\beta+\frac{2\lambda}{6-3\lambda}\gamma\leq\frac{6-q}{2q};\\
\mathrm{(A14)'}&\;p\in\left(\frac{3}{2},3\right),\;q\in[2,2p'),\;\lambda\in\left[\frac{6}{5},\frac{6p}{5p-3}\right),\\
&\alpha+\frac{(4p-6)q}{(6-q)p}\beta\leq1,\;\alpha+\frac{(4p-6)\lambda}{(6-3\lambda)p}\gamma\leq1,\;\beta+\gamma<\frac{3}{q}+\frac{3}{\lambda}-3;\\
\mathrm{(A15)'}&\;p\in\left(\frac{3}{2},3\right),\;q\in[2,2p'),\;\lambda\in\left[\frac{6p}{5p-3},2\right],\\
&\alpha+\frac{(4p-6)q}{(6-q)p}\beta\leq1,\;\alpha+2\gamma<\frac{3}{p}+\frac{6}{\lambda}-4,\;\beta+\gamma<\frac{3}{q}+\frac{3}{\lambda}-3;\\
\mathrm{(A16)'}&\;p\in\left(\frac{3}{2},3\right),\;q\in[2p',6],\;\lambda\in\left[1,\frac{6}{5}\right),\\
&\alpha+2\beta<\frac{3}{p}+\frac{6}{q}-2,\;\alpha+\frac{(4p-6)\lambda}{(6-3\lambda)p}\gamma\leq1,\;\beta+\frac{2\lambda}{6-3\lambda}\gamma\leq\frac{6-q}{2q};\\	\mathrm{(A17)'}&\;p\in\left(\frac{3}{2},3\right),\;q\in[2p',6],\;\lambda\in\left[\frac{6}{5},\frac{6p}{5p-3}\right),\\
&\alpha+2\beta<\frac{3}{p}+\frac{6}{q}-2,\;\alpha+\frac{(4p-6)\lambda}{(6-3\lambda)p}\gamma\leq1,\;\beta+\gamma<\frac{3}{q}+\frac{3}{\lambda}-3.
\end{align*}
Especially, when $\alpha=\beta=\gamma=0$, $p,q,\lambda$ satisfy one of the following assumptions
\begin{align*}
\mathrm{(B1)'}&\;p\in\left[3,\frac{9}{2}\right],q\in[1,2p'),\lambda\in\left[1,\frac{6p}{5p-3}\right);\\
\mathrm{(B2)'}&\;p\in\left(\frac{3}{2},3\right),\;q\in[1,2),\;\lambda\in[1,2];\\
\mathrm{(B3)'}&\;p\in\left(\frac{3}{2},3\right),q\in[1,6],\lambda\in\left[1,\frac{6}{5}\right);\\
\mathrm{(B4)'}&\;p\in\left[3,\frac{9}{2}\right],\;q\in[1,2),\;\lambda\in\left[\frac{6p}{5p-3},2\right],\text{ with }\;\frac{1}{p}+\frac{2}{\lambda}\geq\frac{4}{3};\\
\mathrm{(B5)'}&\;p\in\left[3,\frac{9}{2}\right],\;q\in[2p',6],\;\lambda\in\left[1,\frac{6}{5}\right),\text{ with }\;\frac{1}{p}+\frac{2}{q}\geq\frac{2}{3};\\
\mathrm{(B6)'}&\;p\in\left[3,\frac{9}{2}\right],\;q\in[2,2p'),\;\lambda\in\left[\frac{6p}{5p-3},2\right],\text{ with }\;\frac{1}{p}+\frac{2}{\lambda}\geq\frac{4}{3},\;\frac{1}{q}+\frac{1}{\lambda}\geq1;\\
\mathrm{(B7)'}&\;p\in\left[3,\frac{9}{2}\right],\;q\in[2p',6],\;\lambda\in\left[\frac{6}{5},\frac{6p}{5p-3}\right),\text{ with }\;\frac{1}{p}+\frac{2}{q}\geq\frac{2}{3},\;\frac{1}{q}+\frac{1}{\lambda}\geq1;\\
\mathrm{(B8)'}&\;p\in\left[3,\frac{9}{2}\right],\;q\in[2p',6],\;\lambda\in\left[\frac{6p}{5p-3},2\right],\text{ with }\; \frac{1}{p}+\frac{2}{q}\geq\frac{2}{3},\;\frac{1}{p}+\frac{2}{\lambda}\geq\frac{4}{3},\\
&\frac{1}{q}+\frac{1}{\lambda}\geq1;\\
\mathrm{(B9)'}&\;p\in\left(\frac{3}{2},3\right),\;q\in[2,2p'),\;\lambda\in\left[\frac{6}{5},2\right],\text{ with }\;\frac{1}{q}+\frac{1}{\lambda}\geq1;\\
\mathrm{(B10)'}&\;p\in\left(\frac{3}{2},3\right),\;q\in[2,6],\;\lambda\in\left[\frac{6}{5},\frac{6p}{5p-3}\right),\text{ with }\;\frac{1}{q}+\frac{1}{\lambda}\geq1.\\
\end{align*}
It is worthy mentioning that the assumptions $\mathrm{(A1)}$-$\mathrm{(A17)}$ and $\mathrm{(A1)'}$-$\mathrm{(A17)'}$, $\mathrm{(B1)}$-$\mathrm{(B10)}$ and $\mathrm{(B1)'}$-$\mathrm{(B10)'}$ have a close connection through the transformation
$$\frac{1}{r}=\frac{1}{\lambda}-\frac{1}{3}.$$

The first lemma is the existence and boundedness of the Bogovskii map, which is used to deal with the pressure term $\nabla \pi$.
\begin{Lem}\label{Lem2.1}$($See \cite[Lemma 1]{Tsai21}$)$
Let $E$ be a bounded Lipschitz domain in $\mathbb{R}^{3}$. Denote $L_{0}^{q}(E):=\{f\in L^{q}(E):\int_E fdx=0\}$ with $1<q<\infty$. There exists a linear operator
\begin{equation*}
\mathrm{Bog}:L_{0}^{q}(E)\rightarrow W_{0}^{1,q}(E),
\end{equation*}
such that for any  $f\in L_{0}^{q}(E),w=\mathrm{Bog}f$ is a vector field satisfying
\begin{equation*}
w\in W_{0}^{1,q}(E), \hspace{0.3cm}\mathrm{div}w=f,\hspace{0.3cm}\|\nabla w\|_{L^{q}(E)}\leq C_{\mathrm{Bog}}(E,q)\|f\|_{L^{q}(E)},
\end{equation*}
where the constant $C_{\mathrm{Bog}}(E,q)$  is independent of $f$. By a rescaling argument, we see
$$C_{\mathrm{Bog}}(RE,q)=C_{\mathrm{Bog}}(E,q),\text{ where $RE=\{Rx:x\in E\}$.}$$
\end{Lem}

The second lemma is the so-called Poincar\'{e}-Sobolev inequality. Let $f_{\Omega}$ represent the mean value of the function $f$ over $\Omega$.
\begin{Lem}\label{Lem2.2}$($See \cite[Theorem 3.15]{Giusti} and \cite[Remark 2]{CIY23}$)$
Let $\Omega\subset\mathbb{R}^n$ be a bounded Lipschitz domain. There exists a positive constant $C(n,p,\Omega)$ such that if $1\leq p<n$, then for every $f\in W^{1,p}(\Omega)$, we have
$$\|f-f_{\Omega}\|_{L^{\frac{np}{n-p}}(\Omega)}\leq C(n,p,\Omega)\|\nabla f\|_{L^{p}(\Omega)}.$$
In particular, when $\Omega=B_{2R}\backslash B_R$, by using a rescaling argument, we see that the constant $C(n,p,B_{2R}\backslash B_R)$ does not depend on $R$, i.e.,  $C(n,p,B_{2R}\backslash B_R)=C(n,p)$.
\end{Lem}
Next, we introduce the following standard iteration lemma, which is a generalization and a direct consequence of \cite[Lemma 3.1]{Giaquinta}, and can be found in \cite[Lemma 2.1]{CL24} or \cite[Lemma 2.2]{YX20}.
\begin{Lem}\label{Lem2.3}
Let $f(t)$ be a non-negative bounded function on $\left[r_0, r_1\right] \subset \mathbb{R}^{+}$. If there are non-negative constants $A_i, B_i,\alpha_i$, $i=1,2,\cdots,m$, and a parameter $\theta_0 \in[0,1)$ such that for any $r_0 \leq s<t \leq r_1$, it holds that
$$
f(s) \leq \theta_0 f(t)+\sum_{i=1}^m\left(\frac{A_i}{(t-s)^{\alpha_i}}+B_i\right),
$$
then
$$
f(s) \leq C\sum_{i=1}^m\left(\frac{A_i}{(t-s)^{\alpha_i}}+B_i\right),
$$
where $C$ is a constant depending on $\alpha_1,\alpha_2,\cdots,\alpha_m$ and $\theta_0$.
\end{Lem}

Let $\frac{3R}{2}\leq s<t\leq 2R$, then we infer $R<\frac{2t}{3}<\frac{3t}{4}\leq s$. We introduce a cut-off function $\eta \in C_0^{\infty}(\mathbb{R}^3)$ satisfying
\begin{equation*}
\eta(x)= \begin{cases}1, & |x| <s, \\ 0, & |x| >t,\end{cases}
\end{equation*}
with
$$\text{$0\leq\eta (x)\leq 1$, and $\|\nabla \eta\|_{L^{\infty}} \leq \frac{C}{t-s}$, $\|\nabla^2\eta\|_{L^{\infty}} \leq \frac{C}{(t-s)^2}$.}$$
By Lemma \ref{Lem2.1}, there exists $w\in W_{0}^{1,r}(B_t\backslash B_\frac{2t}{3})$ such that $w$ satisfies the following equation
\begin{equation}\label{ine2.1}
\mathrm{div }w=u\cdot\nabla\eta^{2} \text{ in }B_t\backslash B_\frac{2t}{3},
\end{equation}
with the estimate
\begin{equation}\label{ine2.2}
\aligned
\|\nabla w\|_{L^{r}(B_t\backslash B_\frac{2t}{3})}\leq C\|u\cdot\nabla\eta^2\|_{L^{r}(B_t\backslash B_\frac{2t}{3})}\leq\frac{C}{t-s}\|u\|_{L^{r}(B_t\backslash B_s)},
\endaligned
\end{equation}
for any $1<r<+\infty$. We extend $w$ by zero to $B_\frac{2t}{3}$, then $w\in W_{0}^{1,r}(B_t).$

In order to establish the Liouville type results, we need to estimate several integrals involving $u$, $v$ and $\theta$.
Denote
\begin{align*}
&J_1=\frac{1}{(t-s)^2}\int_{B_t \backslash B_s}(|u|^2+|v|^2+|\theta-\overline{\theta}_{R}|^2) d x,\;J_2=\frac{1}{t-s} \int_{B_t \backslash B_R}|u|^3  d x,\\
&J_3=\int_{B_t\backslash B_\frac{2t}{3}}|v|^{2}|\nabla w|dx+\frac{1}{t-s}\int_{B_t \backslash B_s}|v|^2|u|dx,\\
&J_4=\frac{1}{t-s}\int_{B_t \backslash B_s}|\theta-\overline{\theta}_{R}|^2|u|dx,\;J_5=\frac{1}{t-s}\int_{B_t \backslash B_s}|v||\theta-\overline{\theta}_{R}|dx.
\end{align*}
The estimates on $J_1$, $J_2$, $J_3$, $J_4$, $J_5$ are given in the next five lemmas. The estimates on $J_2$ and $J_3$ come directly from a recent paper \cite {ZZB25} of the second author. The terms $J_1$, $J_4$ and $J_5$ is slightly different from \cite {ZZB25} and \cite{FZ25}.
The only difference is that $\theta$ is replaced by $\theta-\overline{\theta}_{R}$. So the corresponding lemmas for $J_1$, $J_4$ and $J_5$ still hold after the substitution of $\theta$ by $\theta-\overline{\theta}_{R}$.

\begin{Lem}\label{Lem2.4}$($See \cite[Lemma 2.4]{ZZB25}$)$
Let $\frac{3R}{2}\leq s<t\leq 2R$ and $p,q,r\geq1$. Suppose that $u,v,\theta$ are smooth vector-valued functions. Then
\begin{itemize}
\item[(i)] For any $\varepsilon,\delta,\sigma>0$, there exist positive constants $C_\varepsilon$ and $C_\delta$ such that
\begin{equation}\label{ine2.3}
\aligned
J_1\leq&\varepsilon\|u\|_{L^{6}(B_t\backslash B_R)}^{2}+\frac{C_\varepsilon}{(t-s)^{\frac{6}{p}-1}}\|u\|_{L^{p}(B_{2R}\backslash B_R)}^{2}+\frac{C}{(t-s)^2}R^{3-\frac{6}{p}}\|u\|_{L^p (B_{2R}\backslash B_R)}^2+\\
&\delta\|v\|_{L^{6}(B_t\backslash B_R)}^{2}+\frac{C_\delta}{(t-s)^{\frac{6}{q}-1}}\|v\|_{L^{q}(B_{2R}\backslash B_R)}^{2}+\frac{C}{(t-s)^2}R^{3-\frac{6}{q}}\|v\|_{L^{q}(B_{2R}\backslash B_R)}^{2}+\\
&\sigma\|\theta-\overline{\theta}_{R}\|_{L^{6}(B_t\backslash B_R)}^{2}+\frac{C_\sigma}{(t-s)^{\frac{6}{r}-1}}\|\theta-\overline{\theta}_{R}\|_{L^{r}(B_{2R}\backslash B_R)}^{2}\\
&+\frac{C}{(t-s)^2}R^{3-\frac{6}{r}}\|\theta-\overline{\theta}_{R}\|_{L^{r}(B_{2R}\backslash B_R)}^{2}.
\endaligned
\end{equation}
\item[(ii)] It holds that
\begin{equation}\label{ine2.4}
\aligned
J_1\leq\frac{CR^2}{(t-s)^2}\left(\|u\|_{L^6(B_{2R}\backslash B_R)}^2 +\|v\|_{L^6(B_{2R}\backslash B_R)}^2+\|\theta-\overline{\theta}_{R}\|_{L^6(B_{2R}\backslash B_R)}^2\right).
\endaligned
\end{equation}
\end{itemize}
\end{Lem}

\begin{Lem}$\label{Lem2.5}($See \cite[Lemma 2.5]{ZZB25}$)$
Let $\frac{3R}{2}\leq s<t\leq 2R$. Suppose that $u$ is a smooth vector-valued function. Then we have the following conclusions:
\begin{itemize}
\item[(i)]  Assume $p\geq3$. It holds that
\begin{equation}\label{ine2.5}
J_2\leq\frac{C}{t-s}R^{3-\frac{9}{p}}\|u\|_{L^{p}(B_{2R}\backslash B_R)}^{3}.
\end{equation}
\item[(ii)]  Assume $p\in\left[1,3\right)$. It holds that
\begin{equation}\label{ine2.6}
J_2\leq \frac{1}{t-s}\|u\|_{L^p(B_t\backslash B_R)}^{\frac{3p}{6-p}}\|u\|_{L^6(B_t\backslash B_R)}^{\frac{18-6p}{6-p}}.
\end{equation}
\item[(iii)] Assume $p\in\left(\frac{3}{2},3\right)$. For any $\varepsilon>0$, it holds that
\begin{equation}\label{ine2.7}
\aligned
J_2&\leq\varepsilon\|u\|_{L^6(B_t\backslash B_R)}^2+\frac{C_\varepsilon}{(t-s)^\frac{6-p}{2p-3}}\|u\|_{L^p(B_{2R}\backslash B_R)}^\frac{3p}{2p-3}.
\endaligned
\end{equation}
\end{itemize}
\end{Lem}

\begin{Lem}\label{Lem2.6}$($See \cite[Lemma 2.7]{ZZB25}$)$
Let $\frac{3R}{2}\leq s<t\leq 2R$. Suppose that $u,v$ are smooth vector-valued functions. Then
\begin{itemize}
\item[(i)] Let $\frac{3}{2}<p\leq\frac{9}{2}$, $1\leq q<2p'$. For any $\delta>0$, there exists a positive constant $C_\delta$ such that
\begin{equation}\label{ine2.8}
J_3\leq\delta\|v\|_{L^6(B_t\backslash B_R)}^2+\frac{C_\delta}{(t-s)^\frac{(6-q)p'}{(3-p')q}}\|u\|_{L^p(B_{2R}\backslash B_R)}^\frac{(6-q)p'}{(3-p')q}\|v\|_{L^q(B_{2R}\backslash B_R)}^2.
\end{equation}
\item[(ii)] Let $\frac{3}{2}<p\leq\frac{9}{2}$, $1\leq q<2p'$. It holds that
\begin{equation}\label{ine2.9}
\aligned
J_3&\leq\frac{C}{t-s}\|u\|_{L^p(B_t\backslash B_R)}\|v\|_{L^q(B_t\backslash B_R)}^\frac{2(3-p')q}{(6-q)p'}\|v\|_{L^6(B_t\backslash B_R)}^\frac{12p'-6q}{(6-q)p'}.
\endaligned
\end{equation}
\item[(iii)] Let $\frac{3}{2}<p\leq\frac{9}{2}$, $q\geq 2p'$. It holds that
\begin{equation}\label{ine2.10}
\aligned
J_3&\leq\frac{C}{t-s}R^{3-\frac{3}{p}-\frac{6}{q}}\|v\|_{L^q(B_{2R}\backslash B_R)}^{2}\|u\|_{L^p(B_{2R}\backslash B_R)}.
\endaligned
\end{equation}
\end{itemize}
\end{Lem}

\begin{Lem}\label{Lem2.7}$($See \cite[Lemma 2.8]{ZZB25}$)$
Let $\frac{3R}{2}\leq s<t\leq 2R$. Suppose that $u$ is a smooth vector-valued function and $\theta$ is a smooth function. Then
\begin{itemize}
\item[(i)] Let $\frac{3}{2}<p\leq\frac{9}{2}$, $1\leq r<2p'$. For any $\sigma>0$, there exists a positive constant $C_\sigma$ such that
\begin{equation}\label{ine2.11}
J_4\leq\sigma\|\theta-\overline{\theta}_{R}\|_{L^6(B_t\backslash B_R)}^2+\frac{C_\sigma}{(t-s)^\frac{(6-r)p'}{(3-p')r}}\|u\|_{L^p(B_{2R}\backslash B_R)}^\frac{(6-r)p'}{(3-p')r}\|\theta-\overline{\theta}_{R}\|_{L^r(B_{2R}\backslash B_R)}^2.
\end{equation}
\item[(ii)] Let $\frac{3}{2}<p\leq\frac{9}{2}$, $1\leq r<2p'$. It holds that
\begin{equation}\label{ine2.12}
\aligned
J_4&\leq\frac{1}{t-s}\|u\|_{L^p(B_t\backslash B_R)}\|\theta-\overline{\theta}_{R}\|_{L^r(B_t\backslash B_R)}^\frac{2(3-p')r}{(6-r)p'}\|\theta-\overline{\theta}_{R}\|_{L^6(B_t\backslash B_R)}^\frac{12p'-6r}{(6-r)p'}.
\endaligned
\end{equation}
\item[(iii)] Let $\frac{3}{2}<p\leq\frac{9}{2}$, $r\geq 2p'$. It holds that
\begin{equation}\label{ine2.13}
\aligned
J_4&\leq\frac{C}{t-s}R^{3-\frac{3}{p}-\frac{6}{r}}\|\theta-\overline{\theta}_{R}\|_{L^r(B_{2R}\backslash B_R)}^{2}\|u\|_{L^p(B_{2R}\backslash B_R)}.
\endaligned
\end{equation}
\end{itemize}

\end{Lem}

\begin{Lem}\label{Lem2.8}$($See \cite[Lemma 2.9]{FZ25}$)$
Let $\frac{3R}{2}\leq s<t\leq 2R$. Suppose that $v$ is a smooth vector-valued function and $\theta$ is a smooth function. Then
\begin{itemize}
\item[(i)] Let $1\leq q<2$, $1\leq r<2$. For any $\delta,\sigma>0$, there exists a positive constant $C_{\delta,\sigma}$ such that
\begin{equation}\label{ine2.14}
\aligned
J_5\leq&\delta\|v\|_{L^6(B_t\backslash B_R)}^2+\sigma\|\theta-\overline{\theta}_{R}\|_{L^6(B_t\backslash B_R)}^2+\\
&\frac{C_{\delta,\sigma}}{(t-s)^{\frac{(6-q)(6-r)}{q(6-r)+r(6-q)}}}\|v\|_{L^q(B_{2R} \backslash B_R)}^{\frac{2q(6-r)}{q(6-r)+r(6-q)}}\|\theta-\overline{\theta}_{R}\|_{L^r(B_{2R} \backslash B_R)}^{\frac{2r(6-q)}{q(6-r)+r(6-q)}}.
\endaligned
\end{equation}

\item[(ii)] Let $1\leq q<2$, $1\leq r<2$. It holds that
\begin{equation}\label{ine2.15}
J_5\leq\frac{1}{t-s}\|v\|_{L^q(B_t\backslash B_R)}^{\frac{2q}{6-q}}\|\theta-\overline{\theta}_{R}\|_{{L^r}(B_t \backslash B_R)}^{\frac{2r}{6-r}} \|v\|_{L^6(B_t\backslash B_R)}^{\frac{3(2-q)}{6-q}}\|\theta-\overline{\theta}_{R}\|_{{L^6}(B_t \backslash B_R)}^\frac{3(2-r)}{6-r}.
\end{equation}

\item[(iii)] Let $1\leq q<2$, $2\leq r\leq 6$. For any $\delta >0$, there exist a positive constants $C_\delta$ such that
\begin{equation}\label{ine2.16}
\aligned
J_5\leq\delta\|v\|_{L^6(B_t\backslash B_R)}^2+\frac{C_\delta}{(t-s)^\frac{2(6-q)}{6+q}}R^{\frac{3(r-2)(6-q)}{r(6+q)}}\|v\|_{L^q(B_{2R}\backslash B_R)}^{\frac{4q}{6+q}}\|\theta-\overline{\theta}_{R}\|_{L^r(B_{2R} \backslash B_R)}^\frac{2(6-q)}{6+q}.
\endaligned
\end{equation}

\item[(iv)] Let $1\leq q<2$, $2\leq r\leq6$. It holds that
\begin{equation}\label{ine2.17}
J_5\leq\frac{C}{t-s}R^{3(\frac{1}{2}-\frac{1}{r})}\|v\|_{L^q(B_t\backslash B_R)}^{\frac{2q}{6-q}}\|\theta-\overline{\theta}_{R}\|_{{L^r}(B_t \backslash B_R)} \|v\|_{L^6(B_t\backslash B_R)}^{\frac{3(2-q)}{6-q}}.
\end{equation}

\item[(v)] Let $2\leq q\leq6$, $1\leq r<2$. For any $\sigma>0$, there exist a positive constants $C_{\sigma}$ such that
\begin{equation}\label{ine2.18}
\aligned
J_5\leq\sigma\|\theta-\overline{\theta}_{R}\|_{L^6(B_t\backslash B_R)}^2+\frac{C_\sigma}{(t-s)^\frac{2(6-r)}{6+r}}R^{\frac{3(q-2)(6-r)}{q(6+r)}}\|v\|_{L^q(B_{2R}\backslash B_R)}^{\frac{2(6-r)}{6+r}}\|\theta-\overline{\theta}_{R}\|_{L^r(B_{2R} \backslash B_R)}^\frac{4r}{6+r}.
\endaligned
\end{equation}

\item[(vi)] Let $2\leq q\leq6$, $1\leq r<2$. It holds that
\begin{equation}\label{ine2.19}
J_5\leq\frac{C}{t-s}R^{3(\frac{1}{2}-\frac{1}{q})}\|v\|_{L^q(B_t\backslash B_R)}\|\theta-\overline{\theta}_{R}\|_{{L^r}(B_t \backslash B_R)}^{\frac{2r}{6-r}} \|\theta-\overline{\theta}_{R}\|_{L^6(B_t\backslash B_R)}^{\frac{3(2-r)}{6-r}}.
\end{equation}

\item[(vii)] Let $2\leq q\leq6$, $2\leq r\leq 6$. It holds that
\begin{equation}\label{ine2.20}
\aligned
J_5\leq \frac{C}{t-s}R^{3(1-\frac{1}{q}-\frac{1}{r})}\|v\|_{L^q(B_{2R} \backslash B_R)}\|\theta-\overline{\theta}_{R}\|_{L^r(B_{2R} \backslash B_R)}.
\endaligned
\end{equation}
\end{itemize}
\end{Lem}

\section{Proof of Liouville type theorems for the tropical climate model}\label{sec3}

Since all terms involving $\theta$ in \eqref{equ1.1} are the derivatives of $\theta$, the form of \eqref{equ1.1} keeps invariant if we substitute $\theta$ by $\theta-\overline{\theta}_{R}$. Multiply both sides of $\eqref{equ1.1}_{1}$, $\eqref{equ1.1}_{2}$ and $\eqref{equ1.1}_{3}$ by $u \eta^2-w$, $v \eta^2$ and $(\theta-\overline{\theta}_{R})\eta^2$, respectively, integrate over $B_t$ and apply integration by parts. This procedure yields
\begin{equation}\label{ine3.1}
\aligned
&\int_{B_t}\left(|\nabla u|^{2}\eta^2+|\nabla v|^{2}\eta^2+|\nabla \theta|^{2}\eta^2\right) d x\\
=&-\int_{B_t}\left[\nabla u:(u\otimes\nabla \eta^2)+\nabla v:(v\otimes\nabla \eta^2)+\nabla \theta:((\theta-\overline{\theta}_{R})\otimes\nabla \eta^2)\right] d x\\
&+\frac{1}{2} \int_{B_t}(|u|^2+|v|^2+|\theta-\overline{\theta}_{R}|^2) u \cdot \nabla \eta^2 d x+\int_{B_t}\nabla u:\nabla w dx+ \\
&\int_{B_t}(v\otimes v):(u\otimes\nabla \eta^2) dx+\int_{B_t}(\theta-\overline{\theta}_{R}) v\cdot\nabla \eta^2 dx- \int_{B_t}(u\otimes u+v\otimes v): \nabla w  dx.
\endaligned
\end{equation}
Using the Young inequality, we can obtain
\begin{equation}\label{ine3.2}
\aligned
&-\int_{B_t}\left[\nabla u:(u\otimes\nabla \eta^2)+\nabla v:(v\otimes\nabla \eta^2)+\nabla \theta:((\theta-\overline{\theta}_{R})\otimes\nabla \eta^2)\right] d x\\
\leq&\frac{1}{2}\int_{B_t}\left(|\nabla u|^{2}\eta^2+|\nabla v|^{2}\eta^2+|\nabla \theta|^{2}\eta^2\right)dx\\
&+2\int_{B_t}\left(|u\otimes\nabla \eta|^2+|v\otimes\nabla \eta|^2+|(\theta-\overline{\theta}_{R})\otimes\nabla \eta|^2\right)dx.
\endaligned
\end{equation}
Applying the Gagliardo-Nirenberg inequality, we have
\begin{equation}\label{ine3.3}
\aligned
&\|u \eta\|_{L^6 (B_t)}^2+\|v \eta\|_{L^6 (B_t)}^2+\|(\theta-\overline{\theta}_{R})\eta\|_{L^6 (B_t)}^2\\
\leq&C\left(\|\nabla(u \eta)\|_{L^2 (B_t)}^2+\|\nabla(v \eta)\|_{L^2 (B_t)}^2+\|\nabla[(\theta-\overline{\theta}_{R})\eta]\|_{L^2 (B_t)}^2 \right)\\
\leq&C \left(\|\eta\nabla u\|_{L^2 (B_t)}^2+\|\eta\nabla v\|_{L^2 (B_t)}^2+\|\eta\nabla \theta\|_{L^2 (B_t)}^2\right)\\
&+C\left(\|u\otimes\nabla \eta\|_{L^2 (B_t)}^2+\|v\otimes\nabla \eta\|_{L^2 (B_t)}^2+\|(\theta-\overline{\theta}_{R})\otimes\nabla \eta\|_{L^2 (B_t)}^2 \right).
\endaligned
\end{equation}
Combining \eqref{ine3.1}, \eqref{ine3.2} and \eqref{ine3.3}, we have
\begin{equation}\label{ine3.4}
\aligned
&\|u \eta\|_{L^6 (B_t)}^2+\|v \eta\|_{L^6 (B_t)}^2+\|(\theta-\overline{\theta}_{R})\eta\|_{L^6 (B_t)}^2+\|\eta\nabla u\|_{L^2 (B_t)}^2+\|\eta\nabla v\|_{L^2 (B_t)}^2+\|\eta\nabla \theta\|_{L^2 (B_t)}^2\\
&\leq\frac{C}{(t-s)^2}\int_{B_t \backslash B_s}(|u|^2+|v|^2+|\theta-\overline{\theta}_{R}|^2) d x+ \frac{C}{t-s} \int_{B_t \backslash B_s}(|u|^2+|v|^2+|\theta-\overline{\theta}_{R}|^2)|u| d x\\
&+C\int_{B_t\backslash B_\frac{2t}{3}}|\nabla u||\nabla w| dx+C\int_{B_t\backslash B_\frac{2t}{3}}(|u|^{2}+|v|^{2})|\nabla w|dx+\frac{C}{t-s}\int_{B_t \backslash B_s}|v||\theta-\overline{\theta}_{R}|dx.
\endaligned
\end{equation}
As is handled in \cite{ZZB25}, it holds that
\begin{align}\label{ine3.5}
\aligned
\frac{C}{t-s} \int_{B_t\backslash B_s}|u|^3dx+C\int_{B_t\backslash B_{\frac{2t}{3}}}|u|^2|\nabla w| dx
\leq\frac{C}{t-s}\int_{B_t\backslash B_R}|u|^3dx.
\endaligned
\end{align}
Using the Young inequality and \eqref{ine2.2}, we get
\begin{equation}\label{ine3.6}
\aligned
C\int_{B_t\backslash B_\frac{2t}{3}}|\nabla u||\nabla w| dx&\leq\frac{1}{2}\int_{B_t\backslash B_\frac{2t}{3}}|\nabla u|^2dx+C\int_{B_t\backslash B_\frac{2t}{3}}|\nabla w|^2dx\\
&\leq\frac{1}{2}\int_{B_t\backslash B_R}|\nabla u|^2dx+\frac{C}{(t-s)^2}\int_{B_t \backslash B_s}|u|^2dx.
\endaligned
\end{equation}

We denote
\begin{equation*}
f(\rho)=\|u\|_{L^6 (B_\rho)}^2+\|v\|_{L^6 (B_\rho)}^2+\|\theta-\overline{\theta}_{R}\|_{L^6 (B_\rho)}^2+\|\nabla u\|_{L^2 (B_\rho)}^2+\|\nabla v\|_{L^2 (B_\rho)}^2+\|\nabla \theta\|_{L^2 (B_\rho)}^2.
\end{equation*}
Collecting \eqref{ine3.4}, \eqref{ine3.5} and \eqref{ine3.6}, recalling the definitions of the terms $J_1$, $J_2$, $J_3$, $J_4$, $J_5$, we find that
\begin{equation}\label{ine3.7}
\aligned
f(s)\leq&\frac{1}{2}\int_{B_t\backslash B_R}|\nabla u|^2dx+\frac{C}{(t-s)^2}\int_{B_t \backslash B_s}(|u|^2+|v|^2+|\theta-\overline{\theta}_{R}|^2) d x\\
&+ \frac{C}{t-s} \int_{B_t \backslash B_s}(|v|^2+|\theta-\overline{\theta}_{R}|^2)|u| d x+\frac{C}{t-s}\int_{B_t\backslash B_R}|u|^3dx\\
&+C\int_{B_t\backslash B_\frac{2t}{3}}|v|^2|\nabla w|dx+\frac{C}{t-s}\int_{B_t \backslash B_s}|v||\theta-\overline{\theta}_{R}|dx\\
\leq&\frac{1}{2}\int_{B_t\backslash B_R}|\nabla u|^2dx+ CJ_1+CJ_2+CJ_3+CJ_4+CJ_5.
\endaligned
\end{equation}

We divide the assumptions into two main cases, i.e., $3 \leq p\leq\frac{9}{2}$ and $\frac{3}{2} <p<3$.
\begin{proof}[{\bf Proof of Theorem \ref{main1}}]
Let $3 \leq p\leq\frac{9}{2}$. Assume that (A1) holds. Combining \eqref{ine3.7}, \eqref{ine2.3}, \eqref{ine2.5}, \eqref{ine2.8}, \eqref{ine2.11} and  \eqref{ine2.14}, we derive that
\begin{align*}
f(s)\leq&\frac{1}{2}f(t)+\frac{C}{(t-s)^{\frac{6}{p}-1}}\|u\|_{L^p(B_{2R}\backslash B_R)}^{2}+\frac{C}{(t-s)^2}R^{3-\frac{6}{p}}\|u\|_{L^{p}(B_{2R}\backslash B_R)}^{2}\\
&+\frac{C}{(t-s)^{\frac{6}{q}-1}}\|v\|_{L^{q}(B_{2R}\backslash B_R)}^{2}+\frac{C}{(t-s)^2}R^{3-\frac{6}{q}}\|v\|_{L^q(B_{2R}\backslash B_R)}^{2}\\
&+\frac{C}{(t-s)^{\frac{6}{r}-1}}\|\theta-\overline{\theta}_{R}\|_{L^{r}(B_{2R}\backslash B_R)}^{2}+\frac{C}{(t-s)^2}R^{3-\frac{6}{r}}\|\theta-\overline{\theta}_{R}\|_{L^r(B_{2R}\backslash B_R)}^{2}\\
&+\frac{C}{(t-s)}R^{3-\frac{9}{p}}\|u\|_{L^p(B_{2R}\backslash B_R)}^3+\frac{C}{(t-s)^\frac{(6-q)p'}{(3-p')q}}\|u\|_{L^p(B_{2R}\backslash B_R)}^\frac{(6-q)p'}{(3-p')q}\|v\|_{L^q(B_{2R}\backslash B_R)}^2\\
&+\frac{C}{(t-s)^\frac{(6-r)p'}{(3-p')r}}\|u\|_{L^p(B_{2R}\backslash B_R)}^\frac{(6-r)p'}{(3-p')r}\|\theta-\overline{\theta}_{R}\|_{L^r(B_{2R}\backslash B_R)}^2\\
&+\frac{C}{(t-s)^{\frac{(6-q)(6-r)}{q(6-r)+r(6-q)}}}\|v\|_{L^q(B_{2R} \backslash B_R)}^{\frac{2q(6-r)}{q(6-r)+r(6-q)}}\|\theta-\overline{\theta}_{R}\|_{L^r(B_{2R} \backslash B_R)}^{\frac{2r(6-q)}{q(6-r)+r(6-q)}}.
\end{align*}
Applying Lemma \ref{Lem2.3} to the above function inequality, and taking $s=\frac{3}{2}R$ and $t=2R$, we conclude that
\begin{align*}
f(R)\leq&f\left(\frac{3}{2}R\right)\leq CR^{1-\frac{6}{p}}\|u\|_{L^p(B_{2R}\backslash B_R)}^2+CR^{1-\frac{6}{q}}\|v\|_{L^{q}(B_{2R}\backslash B_R)}^{2}\\
&+CR^{1-\frac{6}{r}}\|\theta-\overline{\theta}_{R}\|_{L^{r}(B_{2R}\backslash B_R)}^{2}
+CR^{-\frac{(6-q)p'}{(3-p')q}}\|u\|_{L^p(B_{2R}\backslash B_R)}^\frac{(6-q)p'}{(3-p')q}\|v\|_{L^q(B_{2R}\backslash B_R)}^2\\
&+CR^{2-\frac{9}{p}}\|u\|_{L^p(B_{2R}\backslash B_R)}^3+CR^{-\frac{(6-r)p'}{(3-p')r}}\|u\|_{L^p(B_{2R}\backslash B_R)}^\frac{(6-r)p'}{(3-p')r}\|\theta-\overline{\theta}_{R}\|_{L^r(B_{2R}\backslash B_R)}^2\\
&+CR^{-\frac{(6-q)(6-r)}{q(6-r)+r(6-q)}}\|v\|_{L^q(B_{2R} \backslash B_R)}^{\frac{2q(6-r)}{q(6-r)+r(6-q)}}\|\theta-\overline{\theta}_{R}\|_{L^r(B_{2R} \backslash B_R)}^{\frac{2r(6-q)}{q(6-r)+r(6-q)}}.
\end{align*}
Hence, it holds that
\begin{align*}
f(R)\leq&CR^{1-\frac{6}{p}+2\alpha}X_{\alpha,p,R}^2+CR^{1-\frac{6}{q}+2\beta}Y_{\beta,q,R}^2+CR^{1-\frac{6}{r}+2\gamma}\overline{Z}_{\gamma,r,R}^2+CR^{2-\frac{9}{p}+3\alpha}X_{\alpha,p,R}^3\\
&+CR^{-\frac{(6-q)p'}{(3-p')q}+\frac{(6-q)p'}{(3-p')q}\alpha+2\beta}X_{\alpha,p,R}^\frac{(6-q)p'}{(3-p')q}Y_{\beta,q,R}^2+CR^{-\frac{(6-r)p'}{(3-p')r}+\frac{(6-r)p'}{(3-p')r}\alpha+2\gamma}X_{\alpha,p,R}^\frac{(6-r)p'}{(3-p')r}\overline{Z}_{\gamma,r,R}^2\\
&+CR^{-\frac{(6-q)(6-r)}{q(6-r)+r(6-q)}+\frac{2q(6-r)}{q(6-r)+r(6-q)}\beta+\frac{2r(6-q)}{q(6-r)+r(6-q)}\gamma}Y_{\beta,q,R}^{\frac{2q(6-r)}{q(6-r)+r(6-q)}}\overline{Z}_{\gamma,r,R}^{\frac{2r(6-q)}{q(6-r)+r(6-q)}}.
\end{align*}
Letting $R=R_j\rightarrow+\infty$, thanks to
$$1-\frac{6}{p}+2\alpha<0,\;1-\frac{6}{q}+2\beta\leq0,\;1-\frac{6}{r}+2\gamma\leq0,\;2-\frac{9}{p}+3\alpha\leq0,$$
and
\begin{align*}
&-\frac{(6-q)p'}{(3-p')q}+\frac{(6-q)p'}{(3-p')q}\alpha+2\beta=\frac{(6-q)p}{(2p-3)q}\left(-1+\alpha+\frac{(4p-6)q}{(6-q)p}\beta\right)\leq0,\\
&-\frac{(6-r)p'}{(3-p')r}+\frac{(6-r)p'}{(3-p')r}\alpha+2\gamma=\frac{(6-r)p}{(2p-3)r}\left(-1+\alpha+\frac{(4p-6)r}{(6-r)p}\gamma\right)\leq0,
\end{align*}
\begin{align*}
&-\frac{(6-q)(6-r)}{q(6-r)+r(6-q)}+\frac{2q(6-r)}{q(6-r)+r(6-q)}\beta+\frac{2r(6-q)}{q(6-r)+r(6-q)}\gamma\\
&=\frac{(6-q)(6-r)}{q(6-r)+r(6-q)}\left(-1+\frac{2q}{6-q}\beta+\frac{2r}{6-r}\gamma\right)\leq0,
\end{align*}
we get $u,v \in L^6(\mathbb{R}^3)$ and $\nabla u,\nabla v,\nabla \theta\in L^2(\mathbb{R}^3)$. Furthermore, it holds that
\begin{equation}\label{ine3.8}
\aligned
&\lim_{R\rightarrow+\infty}\|u\|_{L^6(B_{2R}\backslash B_R)}=\lim_{R\rightarrow+\infty}\|v\|_{L^6(B_{2R}\backslash B_R)}=0,\\
&\lim_{R\rightarrow+\infty}\|\nabla u\|_{L^2(B_{2R}\backslash B_R)}=\lim_{R\rightarrow+\infty}\|\nabla v\|_{L^2(B_{2R}\backslash B_R)}=\lim_{R\rightarrow+\infty}\|\nabla \theta\|_{L^2(B_{2R}\backslash B_R)}=0.
\endaligned
\end{equation}
Combining \eqref{ine3.7}, \eqref{ine2.4}, \eqref{ine2.5}, \eqref{ine2.8}, \eqref{ine2.11} and \eqref{ine2.15}, we have
\begin{align*}
f(s)&\leq\frac{CR^2}{(t-s)^2}\left(\|u\|_{L^6(B_{2R}\backslash B_R)}^2 +\|v\|_{L^6(B_{2R}\backslash B_R)}^2 +\|\theta-\overline{\theta}_{R}\|_{L^6(B_{2R}\backslash B_R)}^2\right)\\
&+\frac{C}{t-s}R^{3-\frac{9}{p}}\|u\|_{L^{p}(B_{2R}\backslash B_R)}^{3}+\frac{1}{2}\|\nabla u\|_{L^2(B_t\backslash B_R)}^2+\frac{1}{2}\|v\|_{L^6(B_t\backslash B_R)}^2\\
&+\frac{1}{2}\|\theta-\overline{\theta}_{R}\|_{L^6(B_t\backslash B_R)}^2+\frac{C}{(t-s)^\frac{(6-q)p'}{(3-p')q}}\|u\|_{L^p(B_{2R}\backslash B_R)}^\frac{(6-q)p'}{(3-p')q}\|v\|_{L^q(B_{2R}\backslash B_R)}^2\\
&+\frac{C}{(t-s)^\frac{(6-r)p'}{(3-p')r}}\|u\|_{L^p(B_{2R}\backslash B_R)}^\frac{(6-r)p'}{(3-p')r}\|\theta-\overline{\theta}_{R}\|_{L^r(B_{2R}\backslash B_R)}^2\\
&+\frac{C}{t-s}\|v\|_{L^q(B_t\backslash B_R)}^{\frac{2q}{6-q}}\|\theta-\overline{\theta}_{R}\|_{{L^r}(B_t \backslash B_R)}^{\frac{2r}{6-r}} \|v\|_{L^6(B_t\backslash B_R)}^{\frac{3(2-q)}{6-q}}\|\theta-\overline{\theta}_{R}\|_{{L^6}(B_t \backslash B_R)}^\frac{3(2-r)}{6-r}.
\end{align*}
Taking $s=\frac{3}{2}R$ and $t=2R$ in the above inequality, then using the
Poincar{\'{e}}-Sobolev inequality (see Lemma \ref{Lem2.2}), we get
\begin{align*}
f(R)\leq& f\left(\frac{3}{2}R\right)\leq C\left(\|u\|_{L^6(B_{2R}\backslash B_R)}^2+\|v\|_{L^6(B_{2R}\backslash B_R)}^2+\|\nabla\theta\|_{L^2(B_{2R}\backslash B_R)}^2\right)\\
&+\frac{1}{2}\|\nabla u\|_{L^{2}(B_{2R}\backslash B_R)}^2+CR^{2-\frac{9}{p}}\|u\|_{L^p(B_{2R}\backslash B_R)}^3\\
&+CR^{-\frac{(6-q)p'}{(3-p')q}}\|u\|_{L^p(B_{2R}\backslash B_R)}^\frac{(6-q)p'}{(3-p')q}\|v\|_{L^q(B_{2R}\backslash B_R)}^2\\
&+CR^{-\frac{(6-r)p'}{(3-p')r}}\|u\|_{L^p(B_{2R}\backslash B_R)}^\frac{(6-r)p'}{(3-p')r}\|\theta-\overline{\theta}_{R}\|_{L^r(B_{2R}\backslash B_R)}^2\\
&+CR^{-1}\|v\|_{L^q(B_{2R}\backslash B_R)}^{\frac{2q}{6-q}}\|\theta-\overline{\theta}_{R}\|_{{L^r}(B_{2R}\backslash B_R)}^{\frac{2r}{6-r}} \|v\|_{L^6(B_{2R}\backslash B_R)}^{\frac{3(2-q)}{6-q}}\|\nabla\theta\|_{{L^2}(B_{2R}\backslash B_R)}^\frac{3(2-r)}{6-r}.
\end{align*}
Therefore, we have
\begin{align*}
f(R)\leq&C\left(\|u\|_{L^6(B_{2R}\backslash B_R)}^2+\|v\|_{L^6(B_{2R}\backslash B_R)}^2+\|\nabla\theta\|_{L^2(B_{2R}\backslash B_R)}^2\right)\\
&+\frac{1}{2}\|\nabla u\|_{L^{2}(B_{2R}\backslash B_R)}^2+CR^{2-\frac{9}{p}+3\alpha}X_{\alpha,p,R}^3\\
&+CR^{-\frac{(6-q)p'}{(3-p')q}+\frac{(6-q)p'}{(3-p')q}\alpha+2\beta}X_{\alpha,p,R}^\frac{(6-q)p'}{(3-p')q}Y_{\beta,q,R}^2+CR^{-\frac{(6-r)p'}{(3-p')r}+\frac{(6-r)p'}{(3-p')r}\alpha+2\gamma}X_{\alpha,p,R}^\frac{(6-r)p'}{(3-p')r}\overline{Z}_{\gamma,r,R}^2\\
&+CR^{-1+\frac{2q}{6-q}\beta+\frac{2r}{6-r}\gamma}Y_{\beta,q,R}^{\frac{2q}{6-q}}\overline{Z}_{\gamma,r,R}^{\frac{2r}{6-r}} \|v\|_{L^6(B_{2R}\backslash B_R)}^{\frac{3(2-q)}{6-q}}\|\nabla\theta\|_{{L^2}(B_{2R}\backslash B_R)}^\frac{3(2-r)}{6-r}.
\end{align*}
Letting $R=R_j\rightarrow+\infty$ and thanks to \eqref{ine3.8}, we obtain that $u=v=0$ and $\nabla\theta=0$. Thus, $\theta$ is a constant.

In light of the proof of \cite[Theorem 1.1]{FZ25}, the cases of $\mathrm{(A2)}$-$\mathrm{(A9)}$ can be proved in a similar way as the case of $\mathrm{(A1)}$, so we omit the details.

\end{proof}

\begin{proof}[{\bf Proof of Theorem \ref{main2}}]
	Let $\frac{3}{2} <p<3$. Since
	$$\liminf\limits_{R\rightarrow+\infty}\left(X_{\alpha,p,R}+Y_{\beta,q,R}+\overline{Z}_{\gamma,r,R}\right)<+\infty,$$
	there exists a sequence $R_j\nearrow+\infty$ such that
	$$
	\lim\limits_{j\rightarrow+\infty}X_{\alpha,p,R_j}<+\infty,\;\lim\limits_{j\rightarrow+\infty}Y_{\beta,q,R_j}<+\infty,\;\lim\limits_{j\rightarrow+\infty}\overline{Z}_{\gamma,r,R_j}<+\infty.
	$$
Assume that (A10) holds. Combining \eqref{ine3.7}, \eqref{ine2.3}, \eqref{ine2.7}, \eqref{ine2.8}, \eqref{ine2.11} and  \eqref{ine2.14}, we derive that
	\begin{align*}
		f(s)\leq&\frac{1}{2}f(t)+\frac{C}{(t-s)^{\frac{6}{p}-1}}\|u\|_{L^p(B_{2R}\backslash B_R)}^{2}+\frac{C}{(t-s)^2}R^{3-\frac{6}{p}}\|u\|_{L^{p}(B_{2R}\backslash B_R)}^{2}\\
		&+\frac{C}{(t-s)^{\frac{6}{q}-1}}\|v\|_{L^{q}(B_{2R}\backslash B_R)}^{2}+\frac{C}{(t-s)^2}R^{3-\frac{6}{q}}\|v\|_{L^q(B_{2R}\backslash B_R)}^{2}\\
		&+\frac{C}{(t-s)^{\frac{6}{r}-1}}\|\theta-\overline{\theta}_{R}\|_{L^{r}(B_{2R}\backslash B_R)}^{2}+\frac{C}{(t-s)^2}R^{3-\frac{6}{r}}\|\theta-\overline{\theta}_{R}\|_{L^r(B_{2R}\backslash B_R)}^{2}\\
		&+\frac{C}{(t-s)^\frac{6-p}{2p-3}}\|u\|_{L^p(B_{2R}\backslash B_R)}^\frac{3p}{2p-3}+\frac{C}{(t-s)^\frac{(6-q)p'}{(3-p')q}}\|u\|_{L^p(B_{2R}\backslash B_R)}^\frac{(6-q)p'}{(3-p')q}\|v\|_{L^q(B_{2R}\backslash B_R)}^2\\
		&+\frac{C}{(t-s)^\frac{(6-r)p'}{(3-p')r}}\|u\|_{L^p(B_{2R}\backslash B_R)}^\frac{(6-r)p'}{(3-p')r}\|\theta-\overline{\theta}_{R}\|_{L^r(B_{2R}\backslash B_R)}^2\\
		&+\frac{C}{(t-s)^{\frac{(6-q)(6-r)}{q(6-r)+r(6-q)}}}\|v\|_{L^q(B_{2R} \backslash B_R)}^{\frac{2q(6-r)}{q(6-r)+r(6-q)}}\|\theta-\overline{\theta}_{R}\|_{L^r(B_{2R} \backslash B_R)}^{\frac{2r(6-q)}{q(6-r)+r(6-q)}}.
	\end{align*}
Applying Lemma \ref{Lem2.3} to the above function inequality, we obtain
	\begin{align*}
		f(R)\leq&CR^{1-\frac{6}{p}}\|u\|_{L^p(B_{2R}\backslash B_R)}^2+CR^{1-\frac{6}{q}}\|v\|_{L^{q}(B_{2R}\backslash B_R)}^{2}+CR^{1-\frac{6}{r}}\|\theta-\overline{\theta}_{R}\|_{L^{r}(B_{2R}\backslash B_R)}^{2}\\
		&+CR^\frac{p-6}{2p-3}\|u\|_{L^p(B_{2R}\backslash B_R)}^\frac{3p}{2p-3}+CR^{-\frac{(6-q)p'}{(3-p')q}}\|u\|_{L^p(B_{2R}\backslash B_R)}^\frac{(6-q)p'}{(3-p')q}\|v\|_{L^q(B_{2R}\backslash B_R)}^2\\
		&+CR^{-\frac{(6-r)p'}{(3-p')r}}\|u\|_{L^p(B_{2R}\backslash B_R)}^\frac{(6-r)p'}{(3-p')r}\|\theta-\overline{\theta}_{R}\|_{L^r(B_{2R}\backslash B_R)}^2\\
		&+CR^{-\frac{(6-q)(6-r)}{q(6-r)+r(6-q)}}\|v\|_{L^q(B_{2R} \backslash B_R)}^{\frac{2q(6-r)}{q(6-r)+r(6-q)}}\|\theta-\overline{\theta}_{R}\|_{L^r(B_{2R} \backslash B_R)}^{\frac{2r(6-q)}{q(6-r)+r(6-q)}}.
	\end{align*}
	Hence, we deduce that
	\begin{align*}
		f(R)\leq&CR^{1-\frac{6}{p}+2\alpha}X_{\alpha,p,R}^2+CR^{1-\frac{6}{q}+2\beta}Y_{\beta,q,R}^2+CR^{1-\frac{6}{r}+2\gamma}\overline{Z}_{\gamma,r,R}^2+CR^{\frac{p-6}{2p-3}+\frac{3p}{2p-3}\alpha}X_{\alpha,p,R}^\frac{3p}{2p-3}\\
		&+CR^{-\frac{(6-q)p'}{(3-p')q}+\frac{(6-q)p'}{(3-p')q}\alpha+2\beta}X_{\alpha,p,R}^\frac{(6-q)p'}{(3-p')q}Y_{\beta,q,R}^2+CR^{-\frac{(6-r)p'}{(3-p')r}+\frac{(6-r)p'}{(3-p')r}\alpha+2\gamma}X_{\alpha,p,R}^\frac{(6-r)p'}{(3-p')r}\overline{Z}_{\gamma,r,R}^2\\
		&+CR^{-\frac{(6-q)(6-r)}{q(6-r)+r(6-q)}+\frac{2q(6-r)}{q(6-r)+r(6-q)}\beta+\frac{2r(6-q)}{q(6-r)+r(6-q)}\gamma}Y_{\beta,q,R}^{\frac{2q(6-r)}{q(6-r)+r(6-q)}}\overline{Z}_{\gamma,r,R}^{\frac{2r(6-q)}{q(6-r)+r(6-q)}}.
	\end{align*}
Letting $R=R_j\rightarrow+\infty$, we get $u,v\in L^6(\mathbb{R}^3)$ and $\nabla u,\nabla v,\nabla \theta\in L^2(\mathbb{R}^3)$. Furthermore, \eqref{ine3.8} holds.
	Combining \eqref{ine3.7}, \eqref{ine2.4}, \eqref{ine2.6}, \eqref{ine2.9}, \eqref{ine2.12}, \eqref{ine2.15} and the Poincar{\'{e}}-Sobolev inequality, we have
	\begin{align*}
		f(R)\leq&C\left(\|u\|_{L^6(B_{2R}\backslash B_R)}^2 +\|v\|_{L^6(B_{2R}\backslash B_R)}^2+\|\nabla\theta\|_{L^2(B_{2R}\backslash B_R)}^2\right)+\frac{1}{2}\|\nabla u\|_{L^{2}(B_{2R}\backslash B_R)}^2\\
		&+\frac{C}{R}\|u\|_{L^p(B_{2R}\backslash B_R)}^{\frac{3p}{6-p}}\|u\|_{L^6(B_{2R}\backslash B_R)}^{\frac{18-6p}{6-p}}+\frac{C}{R}\|u\|_{L^p(B_{2R}\backslash B_R)}\|v\|_{L^q(B_{2R}\backslash B_R)}^\frac{2(3-p')q}{(6-q)p'}\|v\|_{L^6(B_{2R}\backslash B_R)}^\frac{12p'-6q}{(6-q)p'}\\
		&+\frac{C}{R}\|u\|_{L^p(B_{2R}\backslash B_R)}\|\theta-\overline{\theta}_{R}\|_{L^r(B_{2R}\backslash B_R)}^\frac{2(3-p')r}{(6-r)p'}\|\nabla\theta\|_{L^2(B_{2R}\backslash B_R)}^\frac{12p'-6r}{(6-r)p'}\\
		&+\frac{C}{R}\|v\|_{L^q(B_{2R}\backslash B_R)}^{\frac{2q}{6-q}}\|\theta-\overline{\theta}_{R}\|_{{L^r}(B_{2R} \backslash B_R)}^{\frac{2r}{6-r}} \|v\|_{L^6(B_{2R}\backslash B_R)}^{\frac{3(2-q)}{6-q}}\|\nabla\theta\|_{{L^2}(B_{2R}\backslash B_R)}^\frac{3(2-r)}{6-r}.
	\end{align*}
	Consequently, we derive that
	\begin{align*}
		f(R)\leq&C\left(\|u\|_{L^6(B_{2R}\backslash B_R)}^2 +\|v\|_{L^6(B_{2R}\backslash B_R)}^2+\|\nabla\theta\|_{L^2(B_{2R}\backslash B_R)}^2\right)+\frac{1}{2}\|\nabla u\|_{L^{2}(B_{2R}\backslash B_R)}^2\\
		&+CR^{\frac{3p}{6-p}\alpha-1}X_{\alpha,p,R}^{\frac{3p}{6-p}}\|u\|_{L^6(B_{2R}\backslash B_R)}^{\frac{18-6p}{6-p}}+CR^{-1+\alpha+\frac{2(3-p')q}{(6-q)p'}\beta}X_{\alpha,p,R}Y_{\beta,q,R}^\frac{2(3-p')q}{(6-q)p'}\|v\|_{L^6(B_{2R}\backslash B_R)}^\frac{12p'-6q}{(6-q)p'}\\
		&+CR^{-1+\alpha+\frac{2(3-p')r}{(6-r)p'}\gamma}X_{\alpha,p,R}\overline{Z}_{\gamma,r,R}^\frac{2(3-p')r}{(6-r)p'}\|\nabla\theta\|_{L^2(B_{2R}\backslash B_R)}^\frac{12p'-6r}{(6-r)p'}\\
		&+CR^{-1+\frac{2q}{6-q}\beta+\frac{2r}{6-r}\gamma}Y_{\beta,q,R}^{\frac{2q}{6-q}}\overline{Z}_{\gamma,r,R}^{\frac{2r}{6-r}} \|v\|_{L^6(B_{2R}\backslash B_R)}^{\frac{3(2-q)}{6-q}}\|\nabla\theta\|_{{L^2}(B_{2R} \backslash B_R)}^\frac{3(2-r)}{6-r}.
	\end{align*}
Letting $R=R_j\rightarrow+\infty$ and using \eqref{ine3.8}, we conclude that $u=v=0$ and $\theta$ is a constant.

In light of the proof of \cite[Theorem 1.2]{FZ25}, the cases of $\mathrm{(A11)}$-$\mathrm{(A17)}$ can be proved in a similar way as the case of $\mathrm{(A10)}$, so we omit the details.
\end{proof}

\begin{Rem}\label{Rem3.1}
Similarly to \cite[Remark 3.1]{FZ25}, we can easily check that
\begin{itemize}
\item[(i)]
The inequality $\beta+\gamma<\frac{3}{q}+\frac{3}{\lambda}-3$ in $\mathrm{(A5)'}$, $\mathrm{(A6)'}$, $\mathrm{(A8)'}$ and $\mathrm{(A9)'}$ can be replaced by
the equality $\beta+\gamma=\frac{3}{q}+\frac{3}{\lambda}-3$, but the price is that we need to assume in addition that
$$\lim\limits_{j\rightarrow+\infty}Y_{\beta,q,R_j}=0\;\text{ or }\lim\limits_{j\rightarrow+\infty}Z'_{\gamma,\lambda,R_j}=0.$$
\item[(ii)]
The inequality $\alpha+2\beta<\frac{3}{p}+\frac{6}{q}-2$ in $\mathrm{(A16)'}$ and $\mathrm{(A17)'}$ can be replaced by the equality $\alpha+2\beta=\frac{3}{p}+\frac{6}{q}-2,$
but the price is that we need to assume in addition that
$$\liminf\limits_{R\rightarrow+\infty}X_{\alpha,p,R}=0,\;\limsup\limits_{R\rightarrow+\infty}(Y_{\beta,q,R}+Z'_{\gamma,\lambda,R})<+\infty,\text{ or }$$
$$\limsup\limits_{R\rightarrow+\infty}(X_{\alpha,p,R}+Z'_{\gamma,\lambda,R})<+\infty,\;\liminf\limits_{R\rightarrow+\infty}Y_{\beta,q,R}=0.$$
\item[(iii)]
The inequality $\alpha+2\gamma<\frac{3}{p}+\frac{6}{\lambda}-4$ in $\mathrm{(A12)'}$ and $\mathrm{(A15)'}$ can be replaced by the equality $\alpha+2\gamma=\frac{3}{p}+\frac{6}{\lambda}-4,$
but the price is that we need to assume in addition that
$$\liminf\limits_{R\rightarrow+\infty}X_{\alpha,p,R}=0,\;\limsup\limits_{R\rightarrow+\infty}(Y_{\beta,q,R}+Z'_{\gamma,\lambda,R})<+\infty,\text{ or }$$
$$\limsup\limits_{R\rightarrow+\infty}(X_{\alpha,p,R}+Y_{\beta,q,R})<+\infty,\;\liminf\limits_{R\rightarrow+\infty}Z'_{\gamma,\lambda,R}=0.$$
\item[(iv)]
The inequality $\beta+\gamma<\frac{3}{q}+\frac{3}{\lambda}-3$ in $\mathrm{(A14)'}$, $\mathrm{(A15)'}$ and $\mathrm{(A17)'}$ can be replaced by the equality $\beta+\gamma=\frac{3}{q}+\frac{3}{\lambda}-3$,
but the price is that we need to assume in addition that
$$\liminf\limits_{R\rightarrow+\infty}Y_{\beta,q,R}=0,\;\limsup\limits_{R\rightarrow+\infty}(X_{\alpha,p,R}+Z'_{\gamma,\lambda,R})<+\infty,\text{ or }$$
$$\limsup\limits_{R\rightarrow+\infty}(X_{\alpha,p,R}+Y_{\beta,q,R})<+\infty,\;\liminf\limits_{R\rightarrow+\infty}Z'_{\gamma,\lambda,R}=0.$$
\end{itemize}
\end{Rem}

\subsection*{Acknowledgements.}
This work was supported by Science Foundation for the Excellent Youth Scholars of Higher Education of Anhui Province Grant No. 2023AH030073 and Natural Science Foundation of the Higher Education Institutions of Anhui Province Grant No. 2023AH050478.

\end {document}